\documentclass[10pt,twoside]{article}
\usepackage{graphicx}
\usepackage{amsmath}
\usepackage{amssymb}
\usepackage{Latex-document}

\title{\bf The Teaching of Proof }
\author{Deborah Loewenberg Ball\thanks{School of Education,
University of Michigan, Ann Arbor, MI 48109-1259, USA. E-mail: dball@umich.edu} \quad Celia
Hoyles\thanks{Institute of Education, University of London, 20 Bedford Way, London WC1H OAL, UK. E-mail:
choyles@ioe.ac.uk}\\ Hans Niels Jahnke\thanks{ Fachbereich Mathematik und Informatik, Universit?t Essen, 45117
Essen, Germany. E-mail: njahnke@uni-essen.de} \quad Nitsa Movshovitz-Hadar\thanks{Department of Education in
Science and Technology, Technion-Israel Institute of Technology, Haifa 32000, Israel. E-mail:
nitsa@tx.technion.ac.il} \vspace*{-0.5cm}}
\date{\vspace{-8mm}}

\begin{document}
\maketitle

\thispagestyle{first} \setcounter{page}{907}

\begin{abstract}

\vskip 3mm
This panel draws on research of the teaching of
mathematical proof, conducted in five countries at different
levels of schooling. With a shared view of proof as essential to
the teaching and learning of mathematics, the authors present
results of studies that explore the challenges for teachers in
helping students learn to reason in disciplined ways about
mathematical claims.
\end{abstract}

\vskip 12mm

\section{Introduction} \label{section 1}\setzero

\vskip-5mm \hspace{5mm}

Proof is central to mathematics and as such should be a key
component of mathematics education. This emphasis can be
justified not only because proof is at the heart of mathematical
practice, but also because it is an essential tool for promoting
mathematical understanding.

This perspective is not always unanimously accepted by either
mathematicians or educators. There have been challenges to the
status of proof in mathematics itself, including predictions of
the `death of proof'. Moreover, there has been a trend in many
countries away from using proof in the classroom (for a survey see
Hanna \& Jahnke, 1996).

In contrast to this, the authors of the present paper agree that
proof must be central to mathematics teaching at all grades.
Nevertheless, there are lessons to be learned from the debates
over the role of proof. For many pupils, proof is just a ritual
without meaning. This view is reinforced if they are required to
write proofs according to a certain pattern or solely with
symbols. Much mathematics teaching in the early grades focuses on
arithmetic concepts, calculations, and algorithms, and, then, as
they enter secondary school, pupils are suddenly required to
understand and write proofs, mostly in geometry. Substantial
empirical evidence shows that this curricular pattern is true in
many countries.

Needed is a culture of argumentation in the mathematics classroom
from the primary grades up all the way through college. However,
we need to know more about the difficulties pupils encounter when
they are confronted with proof and the challenges faced by
teachers who seek to make argumentation central to the
mathematics classroom. The epistemological difficulties that
confront students in their first steps into proof can be compared
to those faced by scientists in the course of developing a new
theory. At the beginning, definitions do not exist.  It is not
clear what has to be proved and what can be presupposed. These
problems are interdependent, and researchers (like students) find
themselves in danger of circular reasoning. In the infancy of a
theory, a proof may serve more to test the credibility or the
fruitfulness of an assumption than to establish the truth of a
statement. Only later, when the theory has become mature (or the
student has come to feel at home in a domain), can a proof play
its mathematical function of transferring truth from assumptions
to a theorem.

All in all, work is needed in three areas with regard to the
teaching of proof. We need (1) a more refined perception of the
role and function of proof in mathematics, including studies of
the practices of proving in which active mathematicians engage
(epistemological analysis), (2) a deeper understanding of the
gradual processes and complexities involved in learning to prove
(empirical research) and (3) the development, implementation and
evaluation of effective teaching strategies along with carefully
designed learning environments that can foster the development of
the ability to prove in a variety of levels as from the primary
through secondary grades and up to college level (design
research).

We begin in Section 2 with an analysis of what mathematical proof
might involve in the primary grades. Section 3 gives results of a
longitudinal study on the development of proving abilities in
grades 8 and 9. Section 4 is based upon an empirical
investigation of college level teaching and shows how the natural
habit of referring to an example can be used as a leverage into
the teaching of proof, and section 5 discusses the idea of
'physical mathematics' as an environment for the teaching of
proof.

\section{What does it take to (teach to) reason in the primary grades?$^1$}
\label{section 2} \setzero

\footnotetext[1]{Author: Deborah Loewenberg Ball}

\vskip-5mm \hspace{5mm}

Although the teaching and learning of mathematical reasoning has
often been seen as a focus only beginning in secondary school,
calls for improvements in mathematics education in the U.S. have
increasingly emphasized the importance of proof and reasoning
from the earliest grades (NCTM, 2000, p.56). While some may
regard such a focus on reasoning and proof secondary to the main
curricular goals in mathematics at this level, we consider
reasoning to be a basic mathematical skill. Yet what might
`mathematical reasoning' look like with young children, and what
might it take for teachers to systematically develop students'
capacity for such reasoning? These questions form one strand of
our research on the teaching and learning of elementary school
mathematics.

We define `mathematical reasoning' as a set of practices and norms
that are collective, not merely individual or idiosyncratic, and
that are rooted in the discipline (Ball $\&$ Bass, 2000, 2002;
Hoover, in preparation). Mathematical reasoning can serve as an
instrument of inquiry for discovering and exploring new ideas, a
process that we call the {\it reasoning of inquiry}. Mathematical
reasoning also functions centrally in justifying or proving
mathematical claims, a process that we call the {\it reasoning of
justification}.  It is this latter on which we focus here.

The reasoning of justification in mathematics, as we see it,
rests on two foundations. One foundation is an evolving {\it body
of public knowledge} --- the mathematical ideas, procedures,
methods, and terms that have already been defined and established
within a given reasoning community.  This knowledge provides a
point of departure, and is available for public use by members of
that community in constructing mathematical claims and in seeking
to justify those claims to others. For professional
mathematicians, the base of public knowledge might consist of an
axiom system for some mathematical structure simply admitted as
given, plus a body of previously developed and publicly accepted
knowledge derived from those axioms. Hence, the base of public
mathematical knowledge defines the grain size of the logical
steps which require no further warrant, that is acceptable within
a given context.  The second foundation of mathematical reasoning
is {\it mathematical language}---symbols, terms, notation,
definitions, and representations¡ªand rules of logic and syntax
for their meaningful use in formulating claims and the networks of
relationships used to justify them. `Language' is used here to
refer to the entire linguistic infrastructure that supports
mathematical communication with its requirements for precision,
clarity, and economy of expression. Language is essential for
mathematical reasoning and for communicating about mathematical
ideas, claims, explanations, and proofs. Some disagreements stem
from divergent or unreconciled uses of terminology, whereas
others are rooted in substantive and conflicting mathematical
claims (Crumbaugh, 1998; Lampert, 1998). The ability to
distinguish these requires sensitivity to the nature and role of
language in mathematics.

We have been tracing the development of mathematical reasoning in
a class of Grade 3 students (ages 8 and 9) across an entire
school year using detailed and extensive records of the class:
videotapes of the daily lessons, the students' notebooks and
tests, interviews with students, and the teacher's plans and
notes. By comparing the class's work at different points in time,
we are able to discern growth in the students' skills of and
dispositions toward reasoning. We offer two brief examples here.
Early in the school year, the teacher presented the problem, `I
have pennies (one-cent coins), nickels, (five-cent coins), and
dimes (ten-cent coins) in my pocket. Suppose I pull out two
coins, what amounts of money might I have?' The children worked
to find solutions to this problem:  2¡é, 6¡é, 10¡é, 11¡é, 15¡é,
and 20¡é.  The teacher asked the students whether they have found
\underline {all} the solutions to the problem, and how they know.
Some students seemed uncertain about the question. Other students
offered explanations:  `If you keep picking up different coins,
you will keep getting the same answers,' `If you write down the
answers and think about it some more until you have them all.'
The students believed they had found them all, but it was because
they could not find any more. Their empirical reasoning satisfied
them. Moreover, they had neither other ideas nor methods for
building a logical argument which would allow them to prove that
this problem (as worded) had exactly six solutions.  They also
did not have the mathematical disposition to ask themselves about
the completeness of their results when working on a problem with
finitely many solutions.

In contrast, consider an episode four months later.  Based on
their work with simple addition problems, the third graders had
developed conjectures about even and odd numbers (e.g., an odd
number plus an odd number equals an even number). They generated
long lists of examples for each conjecture:  3 + 5 = 8, 9 + 7 =
16, 9 + 9 = 18, and so on.  Two girls, amidst this work, argued
to their classmates:  `You can't prove that Betsy's conjecture
({\it odd + odd = even}) always works. Because, um, $\cdots$
numbers go on and on forever, and that means odd numbers and even
numbers go on forever, so you couldn't prove that all of them
work.'  The other children became agitated and one of them
pointed out that no other conclusion that the class had reached
had met this standard. Pointing to some posted mathematical
ideas, the product of previous work, one girl questioned:  `We
haven't even tried \underline{them} with all the numbers there is,
so why do you say that \underline {those} work? We haven't tried
those with all the numbers that there ever could be.'  And other
children reported that they had found many examples, and this
showed that the conjecture was true.  But some were worried:  One
student pointed out that there are `some numbers you can't even
pronounce and some numbers you don't even know are there.' A day
later, however, challenged by the two girls' claim, the class
arrived at a proof.  Representing an odd number as a number than
can be grouped in twos with one left over, they were able to show
that when you add two odd numbers, the two ones  left over would
form a new group of two, forming an even sum:

\begin{figure}[ht]
\begin{center}
  \rotatebox{180}{\reflectbox{\includegraphics[scale=1]{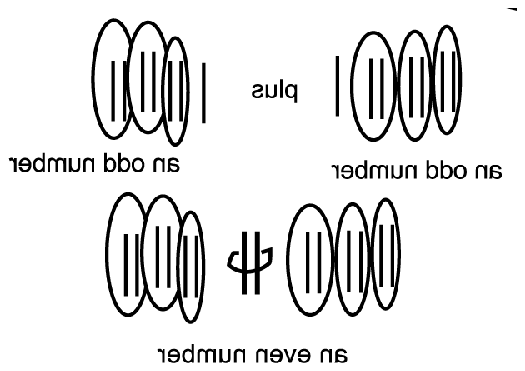}}}

\begin{minipage}{10cm}
\caption{odd + odd = even}
\end{minipage}
\end{center}
\end{figure}

This episode illustrates the important role for definitions.
Having a shared definition for odd and even numbers enabled these
young students to establish a logical argument, based on the
structure of the numbers.  As one girl explained to her
classmates:  `All odd numbers if you circle them by twos, there's
one left over, so if you $\cdots $ plus another odd number, then
the two ones left over will group together, and it will make an
even number.'  The definition equipped them to transcend the
barrier that `numbers go on forever' because it afforded them the
capacity for quantification over an infinite set.  Moreover, this
episode shows the students having developed in their ability to
construct, inspect, and consider arguments using previously
established public mathematical knowledge.

Our research on the nature, structure, and development of
mathematical reasoning has made plain that mathematical reasoning
can be learned, and has highlighted the important role played by
the teacher in developing this capacity. Three domains of work
for the teacher have emerged from our analyses.  A first concerns
the selection of mathematical tasks that create the need and
opportunity for substantial mathematical reasoning.   The
two-coin problem, for instance, did not originally require
students to find all the solutions. Asking this transformed an
ordinary problem into one that involved the need to reason
mathematically about the solution space of the problem. The
second domain of teachers' work centers on making mathematical
knowledge public and in scaffolding the use of mathematical
language and knowledge. Making records of the mathematical work
of the class (through student notebooks, public postings, etc.)
is one avenue, for it helps to make that work public and
available for collective development, scrutiny, and subsequent
use.  This includes attention to where and in what ways knowledge
is recorded, as well as how to name or refer to ideas, methods,
problems, and solutions. Making mathematical knowledge and
language public also requires moving individuals' ideas into the
collective discourse space. A third domain of work, then,
concerns the establishment of a classroom culture permeated with
serious interest in and respect for others' mathematical ideas.
Deliberate attention is required for students to learn to attend
and respond to, as well as use, others' solutions or proposals,
as a means of strengthening their own understanding and the
subsequent contributions they can make to the class's work.

{\bf Acknowledgement.}\  The work reported here draws on my
research with Hyman Bass, Mark Hoover, Jennifer Lewis, and Ed
Wall, as part of the Mathematics Teaching and Learning to Teach
Project at the University of Michigan. This research has been
supported, in part, by grants from the Spencer Foundation and the
National Science Foundation.

\section{The complexity of learning to prove deductively$^2$}\label{section 3} \setzero

\footnotetext[2]{Author: Celia Hoyles}

\vskip-5mm \hspace{5mm}

Deductive mathematical proof offers human beings the purest form
of distinguishing right from wrong; it seems so transparently
straightforward --- yet it is surprisingly difficult for students.
Proof relies on a range of `habits of mind'--- looking for
structures and invariants, identifying assumptions, organising
logical arguments --- each of which, individually, is by no means
trivial.  Additionally these processes have to be coordinated with
visual or empirical evidence and mathematical results and facts,
and are influenced by intuition and belief, by perceptions of
authority and personal conviction, and by the social norms that
regulate what is required to communicate a proof in any particular
situation (see for example, Clements \& Battista, (1992), Hoyles,
(1997), Healy, \& Hoyles, (2000).

The failure of traditional geometry teaching in schools stemmed at
least partly from a lack of recognition of this complexity
underlying proof: the standard practice was simply to present
formal deductive proof (often in a ritualised two-column format)
without regard to its function or how it might connect with
students' intuitions of what might be a convincing argument:
`deductivity was not taught as reinvention, as Socrates did, but
[that it] was imposed on the learner' (Freudenthal, 1973, p.402).
Proving should be part of the problem solving process with
students able to mix deduction and experiment, tinker with ideas,
shift between representations, conduct thought experiments,
sketch and transform diagrams.  But what are the main obstacles
to achieving this flexible habit of mind?

I present here some examples of geometrical questions that have
turned out to be surprisingly difficult--- even for high-
attaining and motivated students.  The analysis forms part of The
Longitudinal Proof Project (Hoyles and K¨¹chemann:
http://www.ioe.ac.uk/proof/), which is analysing students'
learning trajectories in mathematical reasoning over time. Data
are collected through annual surveying of high-attaining students
from randomly selected schools within nine geographically diverse
English regions. Initially 3000 students  (Year 8, age 13) from
63 schools were tested in 2000.  The same students were tested
again in the summer of 2001 using a new test that included some
questions from the previous test together with some new or
slightly modified questions. The same students will be tested
again in June 2002 with the similar aims of testing
understandings and development.

Question G1 in both Year 8 and Year 9 (see Fig 1), is concerned
with how far students use geometrical reasoning to make decisions
in geometry and how far they simply argue from the basis of
perception or what `it looks like' (see Lehrer and Chazan, 1998;
Harel and Sowder, 1998).  In both cases a geometric diagram is
presented, which in the particular case shown, lends support to a
conjecture that turns out to be false. Students are asked whether
or not they agree with the conjecture and to explain their
decision.

Responses to question G1 were coded into 6 broad categories.
Surprisingly, a large number of students in both years simply
answered on the basis of perception and agreed with the false
conjecture with no evidence of progress over the year(Yr 8: 40\%,
Yr 9: 48\%). Additionally 41\% in Yr 8 could come up with a
correct answer and explain this by reference to an explicit
counter-example while only 28\% could do this in Yr 9.

Further analysis, however, is thought-provoking.  Responses to
Yr8G1 showed evidence of three effective strategies: the first to
find the most extreme case that obviously shows that the
diagonals cannot cross at the centre of the circle; the second to
use dynamic reasoning, that is perturbate the diagram in an
incremental way, keeping the given properties invariant (e.g.,
moving one of the vertices round the circumference so the
intersection of the diagonals can no longer be at the centre),
the third is to focus on the diagonals rather than the
quadrilateral and simply to say `I can find opposite vertices such
that the diagonals do not go thru the centre' also evidenced by
students who simply drew a diagonals `cross' without bothering to
draw the quadrilateral itself!  In answer to Yr9G1, it is harder
to find a counter example in a static way as two conditions have
to be controlled (neither diagonal can bisect the area) rather
than only one (one diagonal must not go thru centre); also it is
not possible to find as 'extreme' a counter example as in Y8 (the
nearest equivalent is a concave quadrilateral, though here it is
still possible to end up with two triangles that look very
different but have roughly the same area).  The second strategy
is also harder in the Year 9 question as the dynamic reasoning
has to change an area, not an immediately obvious quantity unlike
the coincidence of two points.  Clearly avoiding the seduction of
perception is only one pitfall in geometrical reasoning.

\begin{center}
  \rotatebox{180}{\reflectbox{\includegraphics[scale=1]{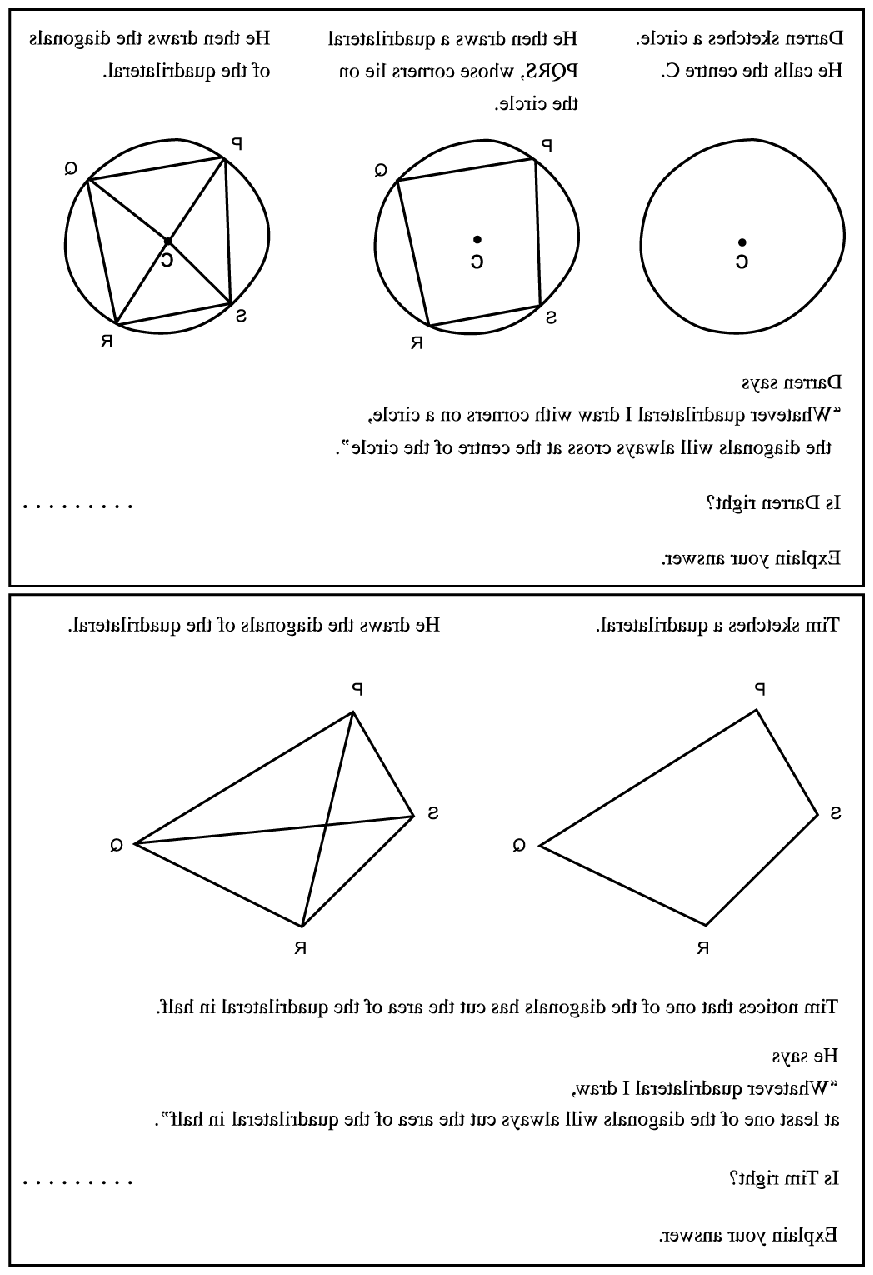}}}

\vskip 2mm

\begin{minipage}{10cm}
Figure 2: Question Y8G1 (top) and Y9G1: Distinguishing perceptual from geometrical reasoning
\end{minipage}
\end{center}

We also found that while students did not easily learn over time
to reject perception, they did improve in calculation. Both Yr 8
and Yr 9 surveys included a question that required knowledge of
certain angle facts (angle on a straight line or at a point,
interior angle sum of triangle; angle property of isosceles
triangle) and where a 3-step calculation had to be performed to
find the size of an angle. We deliberately restricted the task to
working with specific numerical values rather than asking
students to derive a general relationship as would be required in
a standard geometric proof, as this would simply be too hard for
our students who have little experience of proving.

First, we note that students made considerable progress in their
performance onthe calculation part of the question (from 54\%
correin Yr 8 to 73\% in Yr 9). But on analysing responses to the
Yr 9 question, (see Figure 2), where we had asked students to give
reasons for each step of their calculation, we discovered that
not only did they find it hard to match a step in the calculation
to a reason but also they were confused by what it means to give
a reason. Many students interpreted 'reasons' in ways that we did
not anticipate: that is, as an explanation for the step that they
had taken(`u is 40 as I took 40 from 360'), or as request to make
their plans explicit (`I started with p= 320 as the only thing
that I know and I took it from 360 to find u').

Our research is uncovering many more surprises in both student response and progress in proving --- in geometry
but also in algebra (One of our questions (for 14 year-olds) concerns the sum of odd numbers and shows a
remarkably similar spread of responses as those described by Deborah Ball for children age 8/9 years. We know now
even more about potential obstacles to `learning the mathematical game'; but need more systematic work on progress
over time. there are no fool-proof approaches and no short cuts or easy solutions.

{\bf Acknowledgement.}\  I acknowledge the contribution of Dietmar
K¨¹chemann in all the research reported here supported by the
Economic and Social Research Council (ESRC), Project number
R000237777.

\section{The `Because for example $\cdots $'phenomenon, or
transparent pseudo-proofs revisited$^3$} \label{section 4} \setzero

\footnotetext[3]{Author: Nitsa Movshovitz-Hadar}

\vskip-5mm \hspace{5mm}

This panel is about the teaching of proof in mathematics, or as I
interpret it --- providing adequate conditions for gaining
mathematics knowledge. My presentation is based upon the
assumption that mathematics knowledge is {\it in principle} not
different than any other kind of knowledge, although, of course,
the nature of the discipline is different. What, then, is
knowledge? According to Brook and Stainton (2001), a common, long
standing and most plausible answer, given by philosophers to this
question, is that in order to be one's knowledge a proposition
must comply with three necessary (albeit not sufficient)
conditions:
$$
\begin{array}{ll}
(\mbox{i})\quad & \mbox{It must be true.} \\
(\mbox{ii}) \quad &  \mbox{One must believe it. And} \\
(\mbox{iii})\quad& \mbox{One must have justification for believing
it.}
\end{array}
$$

Hanna and Jahnke (1993) suggest that in particular for a novice,
a preliminary step towards appreciating what it is that is being
justified, illumination --- namely understanding and believing, is
of maximum importance. Bertrand Russell makes an important
distinction: Minds do not create truth or falsehood. They create
beliefs. What makes a belief true is its correspondence to a
fact, and this fact does not in any way involve the mind of the
person who holds the belief. This correspondence ensures truth,
and its absence entails falsehood. 'Hence we account
simultaneously for the two facts that beliefs (a) depend on minds
for their {\it existence}, (b) do not depend on minds for their
{\it truth}$\cdots $ ' (Russell, 1912).

We conclude that for a true mathematical statement, i.e., a
theorem, to become one's mathematics knowledge, the learning
environment must consist of teaching tools and strategies that
support the development of two properties: (a) One's belief in
its truth; and (b) one's ability to justify this belief, that is
an ability not just to formally prove it, but also to ensure its
truth by pointing out its correspondence to facts. Said
differently, given a statement p of a mathematical theorem, a
learner should be able to relate to two basic questions: (a) `Do
you believe that p?' and, provided the learner's answer to (a) is
yes, (b)`Why do you believe that p?'

Quite often students' reply to the earlier question is of the
form: `Yes, because for example $\cdots$ '. Very seldom do the
examples that follow, reflect full ability to verify the truth or
even a partial understanding of it.   For example, `Yes, the sum
of every two even integers is an even integer, because for
example 6 plus 8 is 14', does not reflect any insight into the
general case, although it does attest to an understanding of the
statement, (which {\it cannot} be said about the reply: `Yes, for
example, because 14 is the sum of 6 and 8'!) The answer: `Yes,
because for example 6, which is 2$\times$3, plus 8, which is
2$\times$4, give 14, which is 2$\times$7', is slightly better but
not quite. It ties the belief to some acquaintance with the
property of evenness. Although it may be based on deep
understanding, it does {\it not} exhibit more than accepting the
general claim as true, possibly due to a message from an external
authority. (See also Mason 2001, about warrants and the origins
of authority.) To be counted as `satisfactory' the answer should
be something like: `Yes, because for example 6, which is
2$\times$3, plus 8, which is 2$\times$4, give 2$\times$(3+4) and
this IS an even number, as it is a multiple of 2.' This latter one
illustrates what we named {\it a transparent pseudo-proof}.

{\it A transparent proof}, is a proof of a particular case which
is `small enough to serve as a concrete example, yet large enough
to be considered a non-specific representative of the general
case. One can see the general proof through it because nothing
specific to the case enters the proof.' Because a transparent
proof is not a completely polished proof, this kind of `proof' was
later re-named {\it Transparent Pseudo-Proof }or as abbreviated:
{\it Transparent P-Proof}. (Movshovitz-Hadar, 1988, 1998).

The delicate pedagogy involved in preparing a transparent p-proof
was the focus of my ICME-8 Seville presentation (Movshovitz-Hadar
1998). That paper presents the lessons learned through
experimental employment of two slightly different pseudo-proofs,
both of them deserving the title `P-Proof Without Words', yet only
one of which --- `transparent'. The 1998 Samose presentation
(ibid) included further insight into the notion of transparent
p-proof, gained through the preparation of transparent p-proofs as
pedagogical tools to be used in first year linear algebra course,
at Technion - Israel Institute of Technology.

The study of the impact of using transparent p-proofs went on for
four years, and yielded interesting results (Malek, in
preparation). Numerous personal interviews of first-year
mathematics majors and engineering students taking a linear
algebra course, with exposure to transparent p-proofs, yielded
clear evidence as to the impact of reading a transparent p-proof,
on undergraduate students' ability to write, immediately
afterwards, a formal proof of the same claim. A continuing
follow-up also yielded comprehensive evidence as to the impact of
reading transparent p-proofs, on the (passive) ability to read
and comprehend general (formal) proofs, and most important of
all, on the (active) ability to compose general proofs and write
them in a coherent style.

Consequently, we now strongly advocate, wherever it is
appropriate, the use of transparent p-proofs as a pedagogical
tool, as it was shown to support both the development of one's
belief in the truth of mathematical statements and of one's
ability to justify this belief. However, it cannot be
overemphasized that extreme care must be taken by the instructor
in constructing this tool, be it in verbal-symbolic presentation
or in visual-pictorial representation, so that the presentation
is indeed of a transparent proof, namely, it does not hang in any
way to the specifics of the particular case and hence is readily
generalizable. The success of the resulting learning environment
in yielding the development of the ability to prove, depends
heavily on elaborate and careful preparation of the tools by the
instructor.

{\bf Acknowledgement.}\  The research work reported here was
carried by Aliza Malek under my supervision, and was supported by
Technion R \& D funds.

\section{Arguments from physics in mathematical
proofs$^4$}\label{section 5} \setzero

\footnotetext[4]{Author: Hans Niels Jahnke}

\vskip-5mm \hspace{5mm}

Mathematicians often use arguments from physics in mathematical proofs. Some examples, such as the Dirichlet
principle in the variational calculus or \linebreak Archimedes' use of the law of the lever for determining the
volumes of solids, have become famous, and have in fact been regarded by the best mathematicians as elegant
proofs, if not necessarily rigorous. It is only natural, then, that several authors, notably Polya (1954) and
Winter (1978), have proposed that arguments from physics could and should be used in teaching school mathematics.
Besides these publications there are a number of other papers and booklets with examples (see, for example,
Tokieda, 1998). Unfortunately, however, this approach to classroom teaching has not been sufficiently explored.

The application of physics under discussion goes well beyond the
simple physical representation of mathematical concepts, and it
is also distinct from drawing general mathematical conclusions by
the exploration of a large number of instances. Rather, this
approach amounts to using a principle of physics, such as the
uniqueness of the centre of gravity, in a proof and treating it
as if it were an axiom or a theorem of mathematics.

Let us look at a typical example. The so-called Varignon theorem
states that, given an arbitrary quadrangle {\it ABCD}, the
midpoints of its sides W, X, Y, Z form a parallelogram (see
figure 3 below). A purely geometrical proof of this result would
divide the quadrangle into two triangles and apply a similarity
argument.
\begin{figure}[th]
\begin{center}
  \rotatebox{180}{\reflectbox{\includegraphics[scale=1]{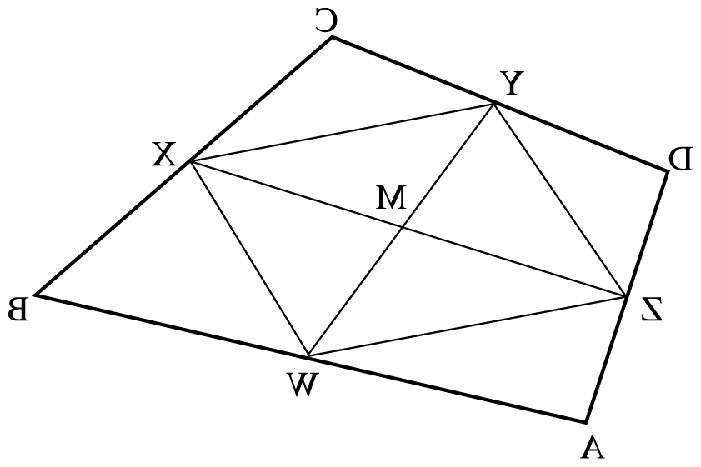}}}

\vskip 2mm

\begin{minipage}{10cm}
\centerline{Figure 3: Varignon's theorem}
\end{minipage}
\end{center}
\end{figure}

An argument from mechanics, on the other hand, would consider
points {\it A, B, C, D} as four weights, each of unit mass,
connected by rigid but weightless rods. Such a system, with a
total mass of 4, has a centre of gravity, and it is this which we
need to determine. The two sub-systems {\it AB} and CD each have
weight 2, and their respective centres of gravity are their
midpoints W and Y. Thus, we may replace {\it AB} and CD by W and Y
loaded with mass 2. Since {\it AB }and CD make up the whole system
{\it ABCD}, its centre of gravity is the midpoint M of WY. In the
same way we can consider {\it ABCD} as made up of {\it BC} and
{\it DA} and its centre of gravity must also be the midpoint of
{\it XZ}. Since the centre of gravity is unique, this midpoint
must be M. This means that M cuts both WY and {\it XZ} into equal
parts. Thus {\it WXYZ}, whose diagonals are WY and {\it XZ}, is a
parallelogram.

The example shows that an argument from physics may
$$
\begin{array}{ll}
o&\quad \hbox{\rm provide a more elegant proof} \\
o &\quad \hbox{\rm reveal the essential features of a complex} \\
& \hskip 3mm \hbox{\rm mathematical structure }\\
 o &\quad \hbox{\rm point out more clearly the relevance of a theorem} \\
 & \hskip 3mm \hbox{\rm to other areas  of mathematics or to other scientific
disciplines}\\
o& \quad \hbox{\rm help create a `holistic' version of a proof,
one that can be grasped in its }\\
&\quad \hbox{\rm entirety,
  as opposed to an elaborate mathematical argument
hard to survey.}
\end{array}
$$

Frequently, an argument from physics helps to {\it generalize} and
to arrive at new theorems. Following the lines of our previous
argument, for example, we can determine the centre of gravity not
only for systems with four masses, but also for those with three,
five, six, and so forth. We can also consider three-dimensional
configurations and investigate whether we are able to translate
the respective statements about the centre of gravity into a
purely geometrical theorem.

There are several reasons why this approach to the teaching of
proof should be further developed and tested. First, it is
unquestionable that, worldwide, we need fresh and possibly more
attractive approaches to the teaching of proof. Since using
arguments from physics in a proof is an alternative to the
established Euclidean routine it might be helpful in motivating
teachers to rethink their attitude to proof.

Another reason is that present-day mathematical practice displays
a significant emphasis on experimentation, and it is only right
that this be reflected in the classroom by a similar emphasis on
experimental mathematics. But it would be dangerous from an
educational point of view if experimental mathematics were to be
represented in the schools only by 'mathematics with computers.'
Quite to the contrary: under the heading of experimental
mathematics, the curriculum should include a strong component
devoted to the classical applications of mathematics to the
physical world. In cultivating this type of mathematics, students
and teachers should be guided by the question of how mathematics
helps to explore and understand the world around us. In this way,
the teaching of proof would be embedded in activities of building
models, inventing arguments to the question 'why', the study of
consequences from assumptions. Working on the border between
mathematics and physics, it could be shown that in quite a few
cases we cannot only apply mathematics to physics, but, vice
versa, can use statements from physics for the derivation of
mathematical theorems.

A Canadian and a German group (Gila Hanna, University of Toronto,
Hans Niels Jahnke, Universit\"at Essen) study the potentials and
pitfalls of this approach in Canadian and German classrooms.
Questions investigated concern the feasibility and the acceptance
of the approach, given the limited knowledge of physics with
students in both countries. It is also asked whether this approach
furthers the general understanding of proof and whether the
students are aware of the difference between using arguments from
physics and the purely empirical appeal to a large number of
instances (Hanna \& Jahnke, 2002).

\label{lastpage}

\end{document}